\documentclass[12pt]{article}
\usepackage{url}
\usepackage[margin=1.1in]{geometry}
\usepackage[noend]{algpseudocode}
\usepackage{algorithm}
\usepackage{times}
\usepackage{amssymb}
\usepackage{amsmath}
\usepackage{latexsym}
\usepackage{graphics}
\usepackage{url}
\usepackage{verbatim}
\usepackage{color}
\usepackage{setspace}
\usepackage{multirow}
\usepackage{ctable}
\usepackage{lineno}
\usepackage{booktabs}
\usepackage{amsfonts}

\algnewcommand\And{\textbf{and}}
\algnewcommand\Or{\textbf{or}}
\algnewcommand\Not{\textbf{not}}
\algnewcommand\In{\textbf{in}}
\algnewcommand\Each{\textbf{each}}
\algnewcommand\Break{\textbf{break}}

\newtheorem{theorem}{Theorem}[section]          

\newcommand{\blue}[1] {\textcolor{blue}{#1}}
\newcommand{\red}[1] {\textcolor{red}{#1}}

\newcommand{\squishlist}{
 \begin{list}{$\bullet$}
  { \setlength{\itemsep}{0pt}
     \setlength{\parsep}{3pt}
     \setlength{\topsep}{3pt}
     \setlength{\partopsep}{0pt}
     \setlength{\leftmargin}{2.5em}
     \setlength{\labelwidth}{1em}
     \setlength{\labelsep}{0.5em} } }

\newcommand{\squishlisttwo}{
 \begin{list}{$\triangleright$}
  { \setlength{\itemsep}{0pt}
     \setlength{\parsep}{0pt}
    \setlength{\topsep}{0pt}
    \setlength{\partopsep}{0pt}
    \setlength{\leftmargin}{2em}
    \setlength{\labelwidth}{1.5em}
    \setlength{\labelsep}{0.5em} } }

\newcommand{\squishend}{
  \end{list}  }

\usepackage{listings}

\definecolor{verbgray}{gray}{0.9}

\lstnewenvironment{code}{%
  \lstset{backgroundcolor=\color{verbgray},
 frame=single,
  language=C,
  framerule=0pt,
  basicstyle=\ttfamily,
  showstringspaces=false,
  commentstyle=\color{blue}\textit,
  keywordstyle=\color{black}\bf,
  numbers=left,
  numberstyle=\footnotesize, 
  columns=fullflexible}}{}

\definecolor{shadecolor}{rgb}{.9, .9, .9}

\newcommand{\inv}{\tau}

\def\Dbar{\leavevmode\lower.6ex\hbox to 0pt
{\hskip-.23ex\accent"16\hss}D}

\def\bZ{{\bf Z}}
\def\bR{{\bf R}}
\def\bC{{\bf C}}
\def\PAF{{\mbox{\rm PAF}}}
\def\PSD{{\mbox{\rm PSD}}}
\def\DFT{{\mbox{\rm DFT}}}

\begin{document}

\title{ Charm bracelets and their application to the construction
of periodic Golay pairs}
\author{
Dragomir {\v{Z}\; \Dbar}okovi{\'c}\thanks{University of Waterloo, Department of Pure Mathematics
and Institute for Quantum Computing, Waterloo, Ontario, N2L 3G1, Canada
e-mail: \texttt{djokovic@uwaterloo.ca}} \ \ \ \
Ilias Kotsireas\thanks{Wilfrid Laurier University,
Department of Physics \& Computer Science,
Waterloo, Ontario, N2L 3C5, Canada
e-mail: \texttt{ikotsire@wlu.ca}} \ \ \ \
Daniel Recoskie\thanks{School of Computer Science, University of
Guelph, Canada. email: \texttt{drecoski@uoguelph.ca}} \ \ \ \
Joe Sawada\thanks{School of Computer Science, University of
Guelph, Canada. Research supported by NSERC. email: \texttt{jsawada@uoguelph.ca}}
}

\date{\today}
\maketitle

\begin{abstract}
A $k$-ary charm bracelet is an equivalence class of length $n$ strings with the action on the indices by the additive group of the ring of integers modulo $n$ extended by the group of units.
By applying an $O(n^3)$ amortized time algorithm to generate charm bracelet representatives with a specified content, we construct 29 new periodic Golay pairs of length $68$.

\end{abstract}

\section{Introduction}

One of the most natural groups acting on $k$-ary strings
$a_0 a_1 \cdots a_{n-1}$ of length $n$ is the group of rotations.
A generator of this group acts on the indices by sending
$i\to i+1 \pmod{n}$, and so sends the string
$a_0 a_1 \cdots a_{n-1} \to a_1 \cdots a_{n-1} a_0$. Applying this action partitions the set of $k$-ary strings into equivalence classes that are called \emph{necklaces}.  When the action of reversal is composed with rotations, the resulting dihedral groups partition $k$-ary strings into equivalence classes called \emph{bracelets}. Generally, we will refer only to the lexicographically smallest element in each respective equivalence class as a necklace or a bracelet.  For example, consider the bracelet equivalence class for the string $12003$:
\begin{center}
\begin{tabular}{cccc}
& 12003  &   30021  \\
& 20031  &   \blue{00213} & $\leftarrow$ \em{bracelet (necklace)} \\
\ \ \ \ \ \ \ \em{necklace} $\rightarrow$ & \blue{00312}  &    02130   \\
& 03120  &    21300  \\
 & 31200  &   13002  \\
\end{tabular}
\end{center}
Observe that this class contains two necklaces 00312 and 00213, the lexicographically smallest being the bracelet representative.

In this paper we generalize the notion of bracelets by considering the action of the group of affine transformations
$j\to a+dj \pmod{n}$ on the indices. Here we consider the indices as elements of the ring of integers modulo $n$ denoted by
$\bZ_n:=\bZ/n\bZ$.  The coefficients $a$ and $d$ also belong to
$\bZ_n$ and $d$ is relatively prime to $n$.   We call the resulting equivalence classes \emph{charm bracelets}. 
Note that if $d \in \{1\}$ we get necklaces, and if $d \in \{1,n-1\}$ we get bracelets.  
As an example, consider the charm bracelet equivalence class for the string $\alpha = a_0a_1a_2a_3a_4 = 12003$:
\begin{center}
\begin{tabular}{cccccc}
 \ \ \ \ \ \ \ \ \ \ \ \ \ \ \ \ \
& 12003 		&10320  & 10230  &   13002  \\
& 20031  		&03201  & 02301  &   30021  \\
& \blue{00312}  &32010  & 23010  &  \blue{00213} &
          $\leftarrow$ \em{charm bracelet}  \\
& 03120         &20103  & 30102  &   02130  \\
& 31200  &\blue{01032}  & \blue{01023}   &  21300   \\
\end{tabular}
\end{center}
Observe that the first strings in each column are the result of the application of the multiplicative group mapping corresponding to $d=1,2,3$ and $4$ respectively.  The subsequent strings in each column correspond to a rotation of the previous string.  Thus, each column will have one necklace representative:
00312, 01032, 01023, and 00213 respectively. The lexicographically smallest necklace 00213 is a charm bracelet. Note that if we take $a=d=n-1$ then the above affine transformation is just the reversal. In general, the maximum number of necklaces in each charm bracelet equivalence class is given by Euler's totient function $\phi(n)$, which denotes the number of positive integers less than $n$ that are relatively prime to $n$. Also, observe that each charm bracelet class will have at most $\phi(n)/2$ bracelets. In particular, observe that the first and last columns of our charm bracelet example correspond to the strings in our previous bracelet example for the string 12003.

Both necklaces and bracelets have been well studied.  Enumeration formulae are well known and efficient algorithms to list necklaces have been given by Fredricksen, Kessler and Maiorana~\cite{fred1,fred2} and Cattell~\emph{et~al}.~\cite{cattell}.  An efficient algorithm to list bracelets is given in \cite{brac_cat}.  Very little is known about charm bracelets except for an  enumeration formula  presented by Titsworth~\cite{titsworth}.  Its binary enumeration sequence was one of the original 2372 sequences presented in 1973 by Sloane in \emph{A Handbook of Integer Sequences}~\cite{handbook}.   In Section~\ref{sec:charm}, we discuss charm bracelets in more detail, presenting a known enumeration formula along with an algorithm to generate them. 

\subsection{An application}

This study of charm bracelets was motivated by the difficult task of deciding the existence of periodic Golay pairs of length 68.  Using our charm bracelet algorithm as step in a searching process we discover 29 new (pairwise nonequivalent) periodic Golay pairs of length 68.  This process is outlined in detail in Section~\ref{sec:app}.

Since our discovery, two separate techniques were discovered to multiply a Golay pair of length $g$ and a periodic Golay pair of length $v$, and obtain as a result a periodic Golay pair of length  $gv$.  We refer loosely to this operation as ``multiplication by $g$''. 
For more details on these multiplications see the recent preprint \cite{PerGol72}.   A special case to multiply by $g=2$  was discovered long ago~\cite[Theorems 2 and 3]{ba90}\footnote{We are grateful to an anonymous referee for pointing this out.}.  Applying the two multiplications by two, the periodic Golay pairs of length 34 presented in~\cite[Theorem 3.1]{FSQ99}) allows us to construct two nonequivalent periodic Golay pair of length 68; however, we have verified that these pairs are not equivalent to any of the $29$ new pairs discovered in this paper (listed in the appendix).  

Finally, we mention that eight non-equivalent periodic Golay pairs of
length 72 have been constructed recently~\cite{PerGol72}.  Consequently, the smallest length for which the existence of periodic Golay pairs is undecided is now 90.

\section{Charm Bracelets}
\label{sec:charm}

\subsection{Enumeration}
An enumeration formula for the number of $k$-ary charm bracelets of length $n$, denoted $CB(n,k)$, was derived in \cite{titsworth}:
\[
CB(n,k) = \frac{1}{n\cdot\phi(n)} \sum_{t=0}^{n-1} \mathop{\sum_{j=1}^{n-1}} \  [\![ \ gcd(n,j)=1  \ ]\!] \ k^{c(j,t)}  \ \mbox{ where }
\]

\[
c(j, t) = \sum_{u=0}^{n-1} \frac{1}{M\left(j, \frac{n}{gcd(n, u(j-1)+t)}\right)}
\]
and where $M(j, L)$ is the smallest positive integer $m$ such that $1+j+\cdots+j^{(m-1)} = 0$ (mod $L$). The Iverson bracket  [\![ $condition$ ]\!] evaluates to 1 if  $condition$ is true, and 0 otherwise.
The enumeration sequence of $CB(n,2)$ corresponds to sequence A002729 in Sloane's {\em The On-Line Encyclopedia of Integer Sequences} \cite{sloane}.  Additionally, the sequences for $CB(n,k)$ for $k=3,4,5$, and 6 correspond to sequences
A056411, A056412, A056413, A056414.

\subsection{Generation algorithm}

Before outlining an algorithm to generate charm bracelets, we first introduce some notation.
Let $\Phi(n)$ denote the set of positive integers less than $n$ that are relatively prime to $n$.
Let $\inv(d,\alpha)$ denote the mapping of $j$ to $dj \bmod n$ acting on the indices of the string $\alpha = a_0a_1\cdots a_{n-1}$.
Let $neck(\alpha)$ denote the necklace representative of the string $\alpha$.
Let $\mathbf{N}_k(n)$ denote the set of all $k$-ary necklaces of length $n$ and let $\mathbf{CB}_k(n)$ denote the set of all $k$-ary charm bracelets of length $n$.

When developing algorithms to exhaustively list combinatorial objects, one of the primary goals is to achieve a CAT algorithm: one that generates each object in constant amortized time.  For charm bracelets this does not appear to be a trivial task.  In this section we outline an algorithm that runs in $O(n^3)$ time per charm bracelet generated.

Perhaps the most straightforward way to exhaustively list $\mathbf{CB}_k(n)$ is by the following approach:
\begin{enumerate}
\item  Generate all the $k$-ary necklaces $\mathbf{N}_k(n)$.

\item  For each necklace $\alpha \in \mathbf{N}_k(n)$  compute  $S(\alpha) = \{ \inv(\alpha, d) \ | \  d \in \Phi(n) \}$.

\item  Compute the necklace of each string  in $S(\alpha)$ to get $T(\alpha) = \{ neck(s) \ | \  s \in S(\alpha) \}$.

\item  Test if $\alpha$ is lexicographically less than or equal to every string in $T(\alpha)$.  If it is, a charm bracelet is found and process $\alpha$.
\end{enumerate}

As mentioned earlier, necklaces can be generated in constant amortized time.  Step 2 requires $O(n^2)$ time to compute the set of $\phi(n)$ strings.  Since the necklace of each string can be computed in $O(n)$ time (see p.222 from ~\cite{combgen}), the set $T$ can also be computed in $O(n^2)$ time.  The third step trivially takes $O(n^2)$ time.  Thus the resulting algorithm runs in $O(n^2)$ time \emph{per necklace}.  Since there are $\phi(n) = O(n)$ necklaces in each charm bracelet class, each charm bracelet gets generated in $O(n^3)$ time.

More detailed pseudocode is given in Algoirthm~\ref{alg:min}.  The function {\sc GenCharm} generates the necklaces using the algorithm from~\cite{cattell,combgen}.  For each necklace $\alpha$ generated, the function {\sc IsCharm}$(\alpha)$ returns whether or not $\alpha$ is a charm bracelet.  It in turn, applies the function {\sc Necklace}$(\beta)$ that returns the necklace of the string $\beta$ by applying a simple modification of the technique given in~\cite{combgen}.    The initial call is {\sc GenCharm}(1,1) initializing $a_0 = 0$.    A complete C implementation is given in the Appendix.

\medskip

\begin{algorithm}[htb]
\small
\caption{Generate all $k$-ary charm bracelets $\alpha = a_1a_2\cdots a_n$ in $O(n^3)$ amortized time. }\label{alg:min}
\begin{algorithmic}[1]
\Statex
\Function{Necklace}{$\beta$}
\State $b_1b_2\cdots b_{2n} \gets \beta \beta$   \ \ \ \   \ \ \ \  \ \ \  \ \blue{$\triangleright$  concatenate two copies of $\beta$}
\State $t \gets j \gets p \gets 1$

\Repeat
       \State $t \gets t+ p  \lfloor \frac{j-t}{p} \rfloor$
	\State $j \gets t+1$
	\State $p \gets 1$
	\While {$j \leq 2n$ {\bf and} $b_{j-p} \le b_j$}
		\If {$b_{j-p} < b_j$} $p  \gets j-t+1$
		\EndIf
		\State $j \gets j + 1$	
	\EndWhile
\Until {$p  \lfloor \frac{j-t}{p} \rfloor  \geq n$}
\State \Return $b_{t} b_{t+1} \cdots b_{t+n-1}$
\EndFunction

\State  \blue{===========================}
\Function{IsCharm}{$\alpha$}
\For {$d \in \Phi(n)$}
	\If { \Call{Necklace}{$ \  \inv(d,\alpha)$ \ } $< \alpha$} \Return {\sc false}
	\EndIf
\EndFor
\State \Return {\sc true}
\EndFunction

\State  \blue{===========================}
\Procedure{GenCharm}{$t,p$}

    \If {$t > n$}
        \If {$N \bmod p = 0$ {\bf and}  \Call{IsCharm}{$\alpha$} } \  \Call{Print}{$\alpha$}  \EndIf
    \Else
        \For {$i$ {\bf from} $a_{t-p}$ {\bf to} $k-1$}
            	\State $a_t \gets i$
            	 \If{$i = a_{t-p}$}  \  \Call{GenCharm}{$t+1,p$}
        	 \Else \  \Call{GenCharm}{$t+1,t$}
       		 \EndIf
        \EndFor
    \EndIf
\EndProcedure
\end{algorithmic}

\end{algorithm}

\newpage

\begin{theorem}
The algorithm {\sc GenCharm} generates all length $n$ charm bracelets in $O(n^3)$-amortized time.
\end{theorem}

As mentioned earlier, the ultimate goal is an algorithm that runs in $O(1)$-amortized time.  However, this appears a very difficult task for charm bracelets.  Any improvement on the $O(n^3)$ algorithm presented here would be a very nice result.   The algorithm can be slightly improved by generating bracelets~\cite{brac_cat} instead of necklaces.   For the application discussed in the next section, only charm bracelets with a specified content are required.  They can also be generated in $O(n^3)$-amortized time by replacing the function {\sc GenCharm}  with the CAT algorithm for fixed content necklaces~\cite{neck_fixed_c} or fixed content bracelets~\cite{brac_fixed}.

\section{Application: Periodic Golay pairs}
\label{sec:app}


Periodic Golay pairs (also known as ``periodic complementary
sequences'') will be defined formally in Section~\ref{sec:pg}.  Early research by Yang
\cite{ya76} used an exhaustive computer search to show that there
are no periodic Golay pairs of length 18. Subsequently, this
case was ruled out by the non-existence result of Arasu and Xiang
\cite{ax92}. For an up-to-date listing of lengths of known
periodic Golay pairs which are not Golay pairs see \cite{PerGol,PerGol72}.  
As mentioned earlier, the smallest length for which the existence of periodic Golay pairs is undecided is now 90.
The periodic Golay pairs can be used to construct Hadamard matrices (see \cite[p. 468]{SY:1992}).

By applying the (fixed-content) charm bracelet algorithm described in the previous section along with a compression of
complementary sequences, we construct 29 periodic Golay pairs of length 68.  One of them will be discussed in more detail
in Section \ref{PerGolExample}. The full listing of the 29 solutions is given in Appendix A.

For the remainder of this section, we use $v$ for the string/sequence lengths rather than the $n$ we used in the previous section, as $v$ is
the standard in design theory.

\subsection{Periodic Golay pairs  vs. Golay pairs}
\label{sec:pg}

The symbols $\bZ,\bR,\bC$ will denote the set of integers,
real numbers and complex numbers, respectively.
Binary sequences will have terms $\pm1$. A pair of binary sequences of length $v$, say,
\begin{equation} \label{seqAB}
A=[a_0,a_1,\ldots,a_{v-1}],\quad
B=[b_0,b_1,\ldots,b_{v-1}]
\end{equation}
is a {\em Golay pair} if for each $k=1,2,\ldots,v-1$:
$$
\sum_{i=0}^{v-k-1} (a_ia_{i+k}+b_ib_{i+k})=0.
$$
It is well known that Golay pairs exist for all lengths
$v=2^a 10^b 26^c$ where $a,b,c$ are nonnegative integers.
For convenience, we shall refer to integers $v$ having this
form as {\em Golay numbers}.
No Golay pairs of other lengths are presently known~
\cite{Borwein:Ferguson:2003}.

We are interested in an analogue of Golay pairs to which we
refer as periodic Golay pairs. They can be defined over any finite abelian group, but we will consider only the finite
cyclic groups. To be specific, we shall use only the cyclic
group $\bZ_v=\{0,1,\ldots,v-1\}$ of integers modulo $v$.
The group operation is addition modulo $v$. From now on we
shall consider the indices of sequences as members of $\bZ_v$.
A {\em periodic Golay pair} is a pair of binary sequences
(\ref{seqAB}) such that for each $k=1,2,\ldots,v-1$:
 \begin{equation} \label{PAF=0}
\sum_{i=0}^{v-1} (a_ia_{i+k}+b_ib_{i+k})=0.
 \end{equation}

Since for any sequence $x_0,x_1,\ldots,x_{v-1}$ we have
$$
\sum_{i=0}^{v-1} x_ix_{i+k} =
\sum_{i=0}^{v-k-1} x_ix_{i+k}+\sum_{i=0}^{k-1} x_ix_{i+v-k},
$$
any Golay pair is also a periodic Golay pair. Therefore
periodic Golay pairs of length $v$ exist whenever $v$ is a
Golay number. However, it is known that they also exist for
some other lengths as well. The first such example was of
length 34 (see \cite{Dj-1998}). At the present time, only finitely many periodic
Golay pairs are known whose length $v$ is not a Golay number.
The smallest length $v$ for which the existence of periodic
Golay pairs of length $v$ is undecided is $v=68$. In this
note we show that such pairs exist.

\subsection{The role of charm bracelets in the search for periodic Golay pairs}

Our objective in this subsection is to explain the role of bracelets in the search
for Golay pairs. In order to do that, we first briefly review some background material.

For an integer sequence $A=[a_0,a_1,\ldots,a_{v-1}]$ of length
$v$, the function $\bZ_v\to\bZ$ which sends
$s\to\sum_{i=0}^{v-1} a_ia_{i+s}$ is known as the {\em periodic
autocorrelation function} (PAF) of $A$. If $(A,B)$ is a
periodic Golay pair of length $v$, then the equation
(\ref{PAF=0}) can be written as
 \begin{equation} \label{PAF=nula}
(\PAF_A+\PAF_B)(s)=0, \quad s=1,2,\ldots,v-1.
 \end{equation}

The {\em discrete Fourier transform} (DFT) of the above
sequence $A$ is the function $\bZ_v\to\bC$ which sends
$s\to\sum_{k=0}^{v-1} a_k \omega^{ks}$, where
$\omega=e^{2\pi i/v}$. The {\em power spectral density} (PSD)
of the sequence $A$ is the function  $\bZ_v\to\bR$
defined by $\PSD_A(s)=|\DFT_A(s)|^2$. By using
\cite[Theorem 2]{Compress}, we deduce that (\ref{PAF=nula})
implies
 \begin{equation} \label{PSD=nula}
(\PSD_A+\PSD_B)(s)=2v, \quad s=0,1,2,\ldots,v-1.
 \end{equation}
 Occasionally we shall write $\PSD(A,s)$ instead of $\PSD_A(s)$,
and similarly for the $\PAF$ function.

Our search for a periodic Golay pair $(A,B)$ is based on the
compression method which is described in detail in the very
recent paper of two of the authors \cite{Compress}. We refer
the reader to this paper also for some additional facts
concerning $A$ and $B$ that we shall use below. In this
computation we used the compression factor $m=2$, and so the
compressed sequences have length $d=v/m=34$. If $a$ and $b$
are the sums of the terms of the sequence $A$ and $B$,
respectively, it is known that $a^2+b^2=4v=136$, and so we may
assume that $a=6$ and $b=10$.

\begin{sloppypar}
In the first stage of the computation we search for suitable
compressed sequences $(A^{(34)},B^{(34)})$. This is a pair of
ternary sequences of length 34,
$$
A^{(34)}=[a_0+a_{34},a_1+a_{35},\ldots,a_{33}+a_{67}], \quad
B^{(34)}=[b_0+b_{34},b_1+b_{35},\ldots,b_{33}+b_{67}],
$$
whose terms $a^{(34)}_i=a_i+a_{i+34}$ and
$b^{(34)}_i=b_i+b_{i+34}$ belong to the set $\{0,2,-2\}$.
Another known fact that we need is that the total number of
0 terms in these two compressed sequences is equal to 34.
For instance, we can choose the case where each of $A$ and $B$
has seventeen 0 terms. As $a=6$ the sequence $A^{(34)}$ must have
the content $(17,10,7)$, i.e., it has seventeen terms equal to 0,
ten terms equal to 2, and seven terms equal to $-2$. Similarly,
$B^{(34)}$ must have the content $(17,11,6)$.
\end{sloppypar}

We can perform on $(A,B)$, as well as on the compressed
sequences, the following operations which preserve the set
of periodic Golay pairs. First, we can permute cyclically
$A$ or $B$ (independently of one another). Second, we can
reverse independently the sequence $A$ or $B$. Third, we can
apply the transformation $x_i\to x_{ki \pmod{v}}$ to both
$A$ and $B$ simultaneously, where $k$ is a fixed integer
relatively prime to $v$. By using these transformations on
the compressed sequences, we deduce that we can restrict
our search for the pairs $(A^{(34)},B^{(34)})$ to the case
where $A^{(34)}$ is a charm bracelet and $B^{(34)}$ is
an ordinary bracelet. (The alphabet used for these bracelets
is $\{0,2,-2\}$.) Since the number of bracelets is much
smaller than the number of all sequences with the same
content, our search will be much faster. There is an
additional speed-up when we restrict (as we may) $A^{(34)}$
to be a charm bracelet. The searches for the bracelets
$A^{(34)}$ and $B^{(34)}$ are performed separately and the
bracelets are written in two files. The search is aborted if
the output file becomes too large. Some of the bracelets do
not need to be recorded. This happens when they fail the
so called PSD test.
In our case this test is based on the fact that we must
have $\PSD(A^{(34)},s)+\PSD(B^{(34)},s)=136$. Hence, the
bracelets for which one of its PSD values is larger than
136 can be safely discarded. By implementing this test
into the search for (charm) bracelets, the size of the output file is considerably reduced.

\subsection{Periodic Golay pairs of length $68$}
\label{PerGolExample}

In this section we present one of the periodic Golay pairs that
we found for length $v = 68$.

Consider the following two sequences of length $34$ each, with $\{-2,0,+2\}$ elements:
\footnotesize
\begin{eqnarray*}
A^{(34)}  & =  & [0,0,0,2,0,0,-2,0,0,0,2,-2,0,0,-2,0,0,2,0,0,0,2,2,-2,0,0,-2,0,0,2,0,2,0,2]  \\
B^{(34)}  & =  & [0,0,-2,2,0,2,0,-2,-2,0,2,2,0,2,-2,0,2,0,-2,2,0,2,2,0,2,0,2,2,0,-2,2,0,-2,-2]
\end{eqnarray*}  \normalsize
These two sequences satisfy the following properties:
\begin{enumerate}
\item $\PAF(A^{(34)},s) + \PAF(B^{(34)},s)=0, s=0,1,\ldots,33$;
\item $\PSD(A^{(34)},s) + \PSD(B^{(34)},s)=2\cdot68=136,~
s = 0,1, \ldots, 33$;
\item $\PSD(A^{(34)},17) = 100$ and $\PSD(B^{(34)},17) = 36$;
\item $\displaystyle\sum_{i=1}^{34} A^{(34)}_i = 6$ and $\displaystyle\sum_{i=1}^{34} B^{(34)}_i = 10$;
\item The total number of $0$ elements in $A^{(34)}$ and $B^{(34)}$ is equal to $34$;
\item The total number of $\pm 2$ elements in $A^{(34)}$ and $B^{(34)}$ is equal to $34$;
\item $A^{(34)}$ contains $21$ zeros and $B^{(34)}$ contains $13$ zeros.
\end{enumerate}

We claim that the sequences $A^{(34)}$ and $B^{(34)}$ are in fact the 2-compressed sequences of two $\{-1,+1\}$ sequences of length $68$ each, that form a particular periodic Golay pair. Here is this particular periodic Golay pair of length $68$:
\begin{eqnarray*}
A &=&
\begin{array}{l}
--++-+-+-++--+--++---++------+-+++ \\
++-++---+-+-+--+-++++++-++-+++++-+ \\
\end{array} \\
B &=&
\begin{array}{l}
---+++---+++++--++-+-+++++++--+--- \\
++-+-++---++-+-++--++++-+-+++-++-- \\
\end{array}
\end{eqnarray*}
In the above periodic Golay pair we use the customary notation of representing $-1$ by $-$ and $+1$ by $+$,
so as to achieve a constant length encoding of the sequences.

In order to find the periodic Golay pair given above, starting from the two sequences $A^{(34)}$ and $B^{(34)}$, we needed to write a program that looks at every individual element of $A^{(34)}$ and $B^{(34)}$ and generates all corresponding potential $\{-1,+1\}$ sequences of length $68$. If we encounter an element equal to $-2$ then this implies that we can set two elements of the length $68$ sequences equal to $-1$. If we encounter an element equal to $+2$ then this implies that we can set two elements of the length $68$ sequences equal to $+1$. If we encounter an element equal to $0$, then this implies that we have two possibilities for the two elements of the length $68$ sequences, either $(-1,+1)$ or $(+1,-1)$. Therefore $A^{(34)}$ generates $2^{21}$ sequences of length $68$ and $A^{(34)}$ generates $2^{13}$ sequences of length $68$. Subsequently we filter these two sets of sequences using the PSD test with PSD constant equal to $136$, since we know from compression theory
\cite{Compress} that the PSD constants of the compressed sequences and the original sequences are equal. The PSD test typically eliminates anywhere between $95$\% to $99$\% of the sequences, so we are left with a very small number of sequences and then it is easy to locate a solution.

Note that there are several thousands (possibly several millions) of pairs of sequences that satisfy 
properties 1 to 6 (and a variant of property 7) of the pair $A^{(34)}, B^{(34)}$, but which do not correspond (via 2-compression) to periodic Golay pairs of order $68$. Both bracelets and charm bracelets are an essential tool for locating such pairs in a systematic manner. On the other hand, all periodic Golay pairs of order $68$ must necessarily be obtained from a pair of sequences of length $34$ that satisfies properties 1 up to 6 and an appropriate version of property 7. Note that property 7 reflects the distribution of the $34$ zeros
in $A^{(34)}, B^{(34)}$ and is directly related with the corresponding bracelets content.

\subsection{Connection with supplementary difference sets}
\label{SDS}

The periodic Golay pairs of fixed length $v$ are in one-to-one correspondence with a special class of combinatorial objects known as supplementary difference sets (SDS). For the definition of SDSs in general we refer the reader to \cite{Compress}. Here
we shall just explain, in the context of this paper, the meaning of SDSs with parameters $(v;r,s;\lambda)=(68;31,29;26)$.
Each of our SDSs consists of two base blocks, say $X$ and $Y$.
They are subsets of the additive group
$\bZ_v=\bZ_{68}=\{0,1,\ldots,67\}$ of sizes $|X|=r=31$ and $|Y|=s=29$. Each nonzero integer in $\bZ_v$ can be represented as
a difference $x_1-x_2$ with $x_1,x_2\in X$ or as a difference
$y_1-y_2$ with $y_1,y_2\in Y$ in total in exactly $\lambda=26$ ways. These particular SDSs are in one-to-one correspondence with
periodic Golay pairs of length $v=68$. Let us make this correspondence explicit. Given an SDS $(X,Y)$ with the above
parameters, we associate to it a periodic Golay pair $(A,B)$. The
first binary sequence $A=[a_0,a_1,\ldots,a_{v-1}]$ is constructed
from the set $X$ by setting $a_j=-1$ if $j\in X$ and $a_j=1$
otherwise. The sequence $B$ is constructed from the set $Y$
in the same way.

We point out that there exist SDSs with two base blocks which do
not correspond to periodic Golay pairs. The SDS's which do correspond to periodic Golay pairs are exactly those whose
parameters satisfy the condition $v=2(r+s-\lambda)$.

\section{Acknowledgements} 
The first two authors wish to acknowledge generous support by NSERC. This work was made possible by the facilities of the Shared Hierarchical Academic Research Computing Network (SHARCNET) and Compute/Calcul Canada. We thank a referee for his suggestions.




\linespread{1.0}
\bibliography{refs}{}
\bibliographystyle{abbrv}

\linespread{1.1}

\newpage
\noindent
\Large
{\bf Appendix A: List of 29 solutions}
\normalsize

\medskip

For convenience we list only the 29 SDSs which correspond to the
29 periodic Golay sequences that we found. All 29 SDSs are given
in the normal form defined in \cite{Dj-2011}. The solution
discussed in Section \ref{PerGolExample} is equivalent to the
solution no. 15 in the list below.

\footnotesize

\begin{eqnarray*}
1) && [[0,1,2,3,4,5,6,7,9,10,12,14,15,20,21,25,28,31,33,34,40,
41,42,45,46,50,52,54,56,57,60], \\
&& [0,2,3,4,6,7,10,11,13,16,18,20,21,23,25,27,28,29,35,36,38,
40,44,45,50,51,58,59,62]],\\
2) && [[0,1,2,3,5,6,7,9,10,11,12,13,17,19,20,21,25,28,31,33,
34,35,40,45,48,49,50,55,58,61], \\
&& [0,1,2,4,7,8,9,12,13,16,18,19,20,22,27,30,35,37,39,41,42,
43,48,50,52,53,56,59,62]],\\
3) && [[0,1,2,3,4,5,6,9,11,15,16,20,22,23,27,29,30,32,36,38,
39,42,43,44,47,48,52,54,55,60,62], \\
&& [0,1,2,3,4,7,8,10,11,12,13,15,16,18,21,22,25,31,33,35,38,
42,44,50,52,55,56,57,60]],\\
4) && [[0,1,2,3,4,5,7,10,11,12,13,15,19,20,21,24,25,27,30,31,
32,37,39,42,46,48,52,55,56,57,59], \\
&& [0,1,2,3,5,6,8,9,10,14,17,20,23,24,27,29,31,33,34,35,39,40,
42,43,47,52,55,57,63]],\\
5) && [[0,1,2,3,4,5,7,11,13,16,19,21,22,27,28,29,30,31,33,35, 38,39,42,43,46,48,49,51,56,58,64], \\
&& [0,1,2,3,7,8,11,12,13,14,15,17,19,21,24,26,27,31,32,35,36,
39,45,46,48,49,53,55,64]],\\
6) && [[0,1,2,3,4,5,8,10,12,13,17,18,19,21,22,24,28,31,32,33, 35,37,38,41,43,45,49,56,57,58,63], \\
&& [0,1,2,3,5,6,7,8,11,12,14,16,20,23,24,26,27,33,34,36,41,42,
45,48,50,52,53,58,65]],\\
7) && [[0,1,2,3,4,5,8,10,14,15,17,23,24,25,26,27,28,29,32,33, 35,36,40,42,43,47,52,54,56,60,62], \\
&& [0,1,2,3,5,7,8,10,13,14,16,18,19,22,25,26,30,31,34,35,36,
39,41,46,49,50,52,56,63]],\\
8) && [[0,1,2,3,4,5,9,11,13,14,15,17,21,22,23,26,28,31,33,35,
38,39,41,42,46,47,50,53,54,56,63],\\
&& [0,1,2,3,4,6,7,9,10,11,12,14,15,20,24,26,27,30, 31,36,41,43,46,47,49,54,56,60,61]],\\
9) && [[0,1,2,3,4,6,7,12,14,15,16,20,22,23,25,26,27,30,32, 34,38,39,40,41,43,44,47,49,52,55,62],\\
&& [0,1,2,4,6,8,9,11,13,14,15,18,19,21,23,29,30,33, 35,36,39,40,44,45,52,53,55,56,63]],\\
10) && [[0,1,2,3,4,6,8,9,10,14,16,17,19,20,23,24,26,28,30,31,35,
36,37,41,45,46,48,49,54,57,58],\\
&& [0,1,2,6,9,10,12,13,15,16,17,19,20,21,25,30,32,35, 37,38,40,42,44,45,51,53,54,56,62]],\\
11) && [[0,1,2,3,4,7,8,10,11,12,13,14,17,21,22,26,27,29,30,33,
38,39,41,46,47,52,54,56,57,58,62],\\
&& [0,1,3,4,5,6,8,9,11,12,15,17,19,21,24,26,27,32,34, 37,38,39,41,46,48,49,53,55,59]],\\
12) && [[0,1,2,3,5,6,7,8,9,14,16,17,20,22,24,26,27,30,31, 33,38,40,43,46,47,48,50,55,58,59,63],\\
&& [0,1,2,4,5,6,7,10,11,15,16,18,19,21,22,27,28,30,34, 36,37,38,40,41,46,48,50,55,61]],\\
13) && [[0,1,2,3,5,6,8,10,12,14,15,17,18,19,23,25,26,29,32, 33,37,40,41,42,43,45,49,52,55,61,62],\\
&& [0,1,2,3,5,6,7,8,10,13,14,18,20,22,23,24,27,28,36, 38,39,41,47,48,49,52,54,58,63]],\\
14) && [[0,1,2,3,5,6,9,10,12,14,15,18,20,21,23,24,25,26,31, 32,33,38,42,43,47,48,50,54,57,58,61],\\
&& [0,1,3,5,6,7,8,9,11,13,15,16,22,25,27,29,30,31,33, 35,40,42,43,44,47,50,56,60,61]],\\
15) && [[0,1,2,3,5,7,8,10,11,17,18,20,21,25,26,27,30,33,34, 38,40,43,44,45,46,47,49,55,57,61,65], \\
&& [0,1,2,3,7,8,10,12,13,14,18,19,21,22,23,25,28,30, 34,36,37,39,40,42,44,47,51,55,56]],\\
\end{eqnarray*}

\begin{eqnarray*}
16) && [[0,1,2,3,5,8,10,11,12,15,19,20,21,24,25,27,28,33,
35,36,39,40,41,43,45,46,47,50,51,57,60],\\
&& [0,1,2,4,5,6,7,10,12,15,18,19,21,22,26,27,28,30, 32,34,36,39,41,43,46,49,55,56,57]],\\
17) && [0,2,4,6,7,8,9,11,12,13,17,18,19,22,23,25,27,33,34, 35,36,37,39,42,43,46,49,51,56,57,59]],\\
&& [[0,1,2,3,7,8,9,11,12,14,15,18,21,25,27,28,29,30, 33,35,39,40,44,47,48,52,55,57,60],\\
18) && [[0,1,2,4,5,8,9,10,11,15,16,19,20,21,22,23,29,31,33, 35,37,40,41,46,48,51,53,55,56,60,63],\\
&& [0,1,2,3,5,6,7,10,11,13,14,16,22,24,26,27,29,30, 32,33,36,39,40,41,45,48,50,56,57]],\\
19) && [[0,1,2,4,5,8,9,11,12,15,17,20,21,22,23,27,28,29,30, 33,37,39,44,45,46,49,50,53,55,58,59],\\
&& [0,1,3,5,6,7,8,9,11,13,14,17,19,20,22,26,27,29,31, 37,40,41,42,44,45,52,54,56,61]],\\
20) && [[0,1,2,5,6,7,8,10,12,13,16,17,18,20,24,28,29,30, 31,33,34,38,40,42,47,48,49,53,59,61,62],\\
&& [0,1,2,4,6,7,9,10,11,13,14,16,19,20,23,26,28,34, 36,37,39,45,49,50,52,53,54,57,60]],\\
21) && [[0,1,2,5,6,7,9,10,11,15,17,19,20,24,26,27,31,32,33, 35,38,39,42,44,45,47,52,55,56,58,59],\\
&& [0,1,2,3,4,5,6,9,10,11,12,13,16,19,22,24,25,27,32, 35,39,40,44,46,50,52,54,56,61]],\\
22) && [0,1,2,3,4,7,8,10,12,13,14,16,19,20,22,23,25,27,32, 33,34,36,38,41,48,49,50,53,54,58,61]],\\
&& [[0,1,3,4,5,6,7,8,11,12,14,19,20,21,26,29,31,35,36, 39,44,45,47,48,50,52,58,62,64],\\
23) && [0,1,2,3,4,6,8,9,10,13,17,21,22,25,26,28,32,33,34, 35,38,40,42,45,48,49,50,51,56,59,62]],\\
&& [[0,1,3,4,5,6,10,13,14,15,18,19,20,22,24,26,29,30,33, 36,38,39,40,41,44,46,51,57,59],\\
24) && [[0,1,3,4,5,7,8,9,10,12,14,15,16,17,24,25,26,27,30, 33,35,36,39,40,43,50,51,55,56,57,63],\\
&& [0,1,2,5,7,8,9,11,13,17,18,20,21,23,25,31,32,35, 36,38,40,42,45,46,49,51,54,57,59]],\\
25) && [[0,1,3,4,5,7,8,11,13,14,15,19,20,21,24,27,29,30, 31,33,34,36,38,39,44,46,47,48,51,57,60],\\
&& [0,2,3,4,6,8,9,10,14,15,17,19,22,24,25,26,28, 35,36,40,43,45,48,49,53,54,55,56,60]],\\
26) && [[0,1,3,4,5,7,9,10,11,12,15,16,18,19,20,26,27,29,31, 33,34,36,38,39,42,43,51,52,56,57,59],\\
&& [0,2,3,4,5,6,10,11,12,13,16,17,22,23,25,27,31, 34,38,40,41,44,46,48,51,54,56,59,60]],\\
27) && [[0,1,3,4,6,7,9,10,11,12,14,18,19,20,22,26,30,32,33, 34,35,37,39,42,47,49,50,51,55,56,60],\\
&& [0,1,2,3,6,9,10,12,14,15,16,17,18,20,23,25,27,30, 34,37,38,43,44,49,50,54,56,59,60]],\\
28) && [[0,2,3,4,5,6,7,9,11,14,15,18,19,21,24,31,32,33,35, 39,40,41,45,47,50,51,52,57,59,60,64],\\
&& [0,1,2,3,4,5,6,9,10,11,15,16,19,22,24,25,26,27,30, 33,36,38,40,44,46,49,53,56,61]],\\
29) && [[0,2,3,4,6,8,9,10,11,12,15,17,18,21,25,28,29,30, 34,38,39,41,44,46,48,49,53,54,56,60,61],\\
&& [0,1,2,3,6,7,8,9,11,13,17,18,21,22,24,27,29,30, 31,33,35,38,41,42,43,52,55,56,58]].
\end{eqnarray*}

\newpage
\noindent
\Large
{\bf Appendix B: C code to generate charm bracelets}
\medskip

\scriptsize
\begin{code}
#include <stdio.h>		
int a[100],b[100];
int N,K,total=0;
//-------------------------------------------------------------
int Gcd(int x, int y){
    int t;

    while( y != 0 ) {
        t = y;  y = x % y;  x = t;
    }
    return x;
}
//-------------------------------------------------------------
void Print() {
    int i;

    total++;
    for (i=1; i<=N; i++) printf("%d", a[i]);
    printf("\n");
}
//-------------------------------------------------------------
// Find the necklace of the string b[1..n] by concatenating two
//   copies of b[1..n] together.  The necklace will be start at
//   index t.  O(n) time.
//-------------------------------------------------------------
int Necklace(){
    int j,t,p;

    for (j=1; j<=N; j++) b[N+j] = b[j];

    j=t=p=1;
    do {
        t = t + p*((j-t)/p);
        j = t + 1;
        p = 1;
        while (j <= 2*N && b[j-p] <= b[j]) {
            if (b[j-p] < b[j]) p = j-t+1;
            j++;
        }
    } while (p * ((j-t)/p) < N);

    return t;
}
//-------------------------------------------------------------
// For each i relatively prime to N, map index j to (ij mod N)
// Then find the necklace of the resulting string, if that
//    necklace is less than the necklace a[1..n] - reject
//-------------------------------------------------------------
int IsCharm(){
    int i,j,offset;

    for(i=2; i<=N-1; i++){
        if ( Gcd(i,N) == 1) {

            // Perform the mapping then determine the necklace
            for(j=0; j<N; j++) b[(j*i)%N + 1] = a[j+1];
            offset = Necklace();

            for (j=1; j<=N; j++){
                if (a[j] < b[offset + j-1]) break;
                else if (a[j] > b[offset + j-1]) return 0;
            }
        }
    }
    return 1;
}
//--------------------------------------------------------------
// Generate necklaces and then check if they are charm bracelets
//--------------------------------------------------------------
int GenCharm(int t, int p) {
    int i;

    if (t > N) {
        if (N%p == 0 && IsCharm()) Print();
    }
    else {
        for (i=a[t-p]; i<K; i++) {
            a[t] = i;
            if (i == a[t-p]) GenCharm(t+1,p);
            else GenCharm(t+1,t);
        }
    }
}
//--------------------------------------------------------------
int main() {

    printf("Enter N K: ");
    scanf("%d %d", &N, &K);

    a[0] = 0;
    GenCharm(1,1);
    printf("Total = %d\n", total);
}
\end{code}


\end{document}